\documentclass[11pt]{amsart}

\usepackage{fullpage}
\usepackage{url}
\usepackage{amssymb}
\usepackage{cite}

\renewcommand{\a}{\alpha}

\newcommand{\g}{\gamma}
\renewcommand{\d}{\delta}
\newcommand{\D}{\Delta}
\newcommand{\e}{\varepsilon}
\newcommand{\f}{\phi}
\newcommand{\s}{\sigma}
\renewcommand{\k}{\kappa}
\renewcommand{\l}{\lambda}

\newcommand{\cC}{{\mathcal C}}
\newcommand{\cM}{{\mathcal M}}
\newcommand{\cN}{{\mathcal N}}
\newcommand{\cB}{{\mathcal B}}
\newcommand{\cK}{{\mathcal K}}

\newcommand{\bR}{\mathbb R}

\newcommand{\bZ}{\mathbb Z}
\newcommand{\bS}{\mathbb S}
\newcommand{\bM}{\mathbb M}
\newcommand{\bH}{\mathbb H}

\newcommand{\Di}{\mathrm{Diff}}

\newcommand{\be}{\begin{equation}}
\newcommand{\ee}{\end{equation}}

\newcommand{\bel}[1]{\begin{equation}\label{#1}}
\newcommand{\beaa}{\begin{eqnarray*}}
\newcommand{\bea}{\begin{eqnarray}}
\newcommand{\beal}[1]{\begin{eqnarray}\label{#1}}
\newcommand{\bean}{\begin{eqnarray}\nonumber}
\newcommand{\beadl}[1]{\begin{deqarr}\label{#1}}
\newcommand{\eeadl}[1]{\arrlabel{#1}\end{deqarr}}
\newcommand{\eeal}[1]{\label{#1}\end{eqnarray}}
\newcommand{\eead}[1]{\end{deqarr}}
\newcommand{\eea}{\end{eqnarray}}
\newcommand{\eeaa}{\end{eqnarray*}}

\newcounter{mnotecount}[section]

\renewcommand{\to}{\rightarrow}

\renewcommand{\exp}{\operatorname{exp}}

\DeclareMathOperator{\Conf}{Conf}
\DeclareMathOperator{\Diff}{Diff}

\renewcommand{\phi}{\varphi}
\renewcommand{\epsilon}{\varepsilon}

\newcommand{\dm}{{\partial M}}

\renewcommand{\hbar}{{\overline h}}

\newcommand{\pre}[2]{{{\vphantom{#2}}^{#1}}\kern-.2ex{#2}}

\theoremstyle{plain}
\newtheorem{theorem}{\sc Theorem}[section]

\newtheorem{lemma}[theorem]{\sc Lemma}

\newtheorem{proposition}[theorem]{\sc Proposition}

\theoremstyle{definition}

\def\endproof{\qed \medskip}
\def\blacksquare{\hbox to .60em {\vrule width .60em height .60em}}

\numberwithin{equation}{section}

\date{\today}

\begin{document}

\title[Boundary value problems on 3-manifolds]{Boundary value problems for metrics on 3-manifolds}

\author[ ]{Michael T. Anderson}
\address{Dept.~of Mathematics, Stony Brook University, Stony Brook, N.Y.~11794-3651, USA} 
\email{anderson@math.sunysb.edu}

\thanks{Partially supported by NSF grant DMS 0905159}

\subjclass{Primary 53C42; Secondary 35J57}
\keywords{conformal immersions, prescribed mean curvature}

\begin{abstract}
We discuss the problem of prescribing the mean curvature and conformal class as boundary 
data for Einstein metrics on 3-manifolds, in the context of natural elliptic boundary value 
problems for Riemannian metrics. 
\end{abstract}

\maketitle

\section{Introduction}\label{section:intro}

  A question long of basic interest to geometers is the existence of 
complete Einstein metrics on manifolds. Any kind of theory for the existence or 
uniqueness of such metrics on compact manifolds is still far from sight. The 
only exception to this is the remarkable work of Perelman and Hamilton, which 
essentially gives a complete theory for closed 3-manifolds. 

  Instead of considering closed manifolds, it might be somewhat simpler to 
consider manifolds with boundary and look for a theory providing existence 
(and uniqueness) for geometrically natural boundary value problems. This 
has recently met with some success, in the context of complete conformally 
compact Einstein metrics, where one prescribes a conformal metric at 
conformal infinity \cite{An1}, and in the context of a natural exterior 
boundary value problem for the static vacuum Einstein equations, \cite{AK}. 

  In this note, we consider the simplest situation, namely boundary value problems 
for Einstein metrics in dimension 3, where the metrics are of constant curvature. 
Seemingly the simplest or most naive question one could ask in this context is 
the following: 

\medskip

{\sc Question}. Given a metric $\g$ on a boundary surface $\dm = S^{2}$ for instance, 
is there an Einstein metric (flat or constant curvature) on the 3-ball $M = B^{3}$ 
inducing $\g$ on $\dm$? 

\medskip

  However, this is basically the isometric immersion problem for surfaces in 
$\bR^{3}$, (or other space-forms), and is a notoriously difficult problem, 
also far from any current resolution. Note however that there are examples 
of smooth metrics on $S^{2}$ which do not isometrically immerse in $\bR^{3}$, 
cf.~\cite{P}, so the answer to the question is no in general.

  The main difficulty here is that although the Einstein equations form an elliptic 
system of equations in a suitable gauge, Dirichlet boundary data for such a 
system never give rise to an elliptic boundary value problem. The Gauss 
constraint equation, (Gauss' Theorema Egregium), is an obstruction to such 
ellipticity. Thus, one should first consider what are the natural boundary 
value problems for the Einstein equations. 

  To describe this, let $M$ be any 3-manifold with boundary $\dm$ which admits a 
metric of constant sectional curvature $\k$. We assume that 
$$\pi_{1}(M, \dm) = 0;$$
by elementary covering space arguments, this means that $\dm$ is connected 
and any loop in $M$ is homotopic to a loop in $\dm$, so that $M$ is a 
3-dimensional handlebody.

  Let $\cM_{\k} = \cM_{\k}^{m,\a}$ be the moduli space of metrics of constant curvature $\k$ 
on $M$ which are $C^{m,\a}$ up to $\dm$, $m \geq 2$, $\a \in (0,1)$. This is the space of 
all such constant curvature metrics $\bM_{\k}$ modulo the action of $\Diff_{1}^{m+1,\a}$ 
of diffeomorphisms of $M$ equal to the identity on $\dm$. In the case of $\k = 0$ for 
instance, the developing map gives an isometric immersion
$$D: (M, g) \to \bR^{3},$$
which induces an isometric Alexandrov immersion of $(\dm, \g)$ into $\bR^{3}$, where $\g = 
g_{\partial M}$. (An immersion of a surface in $\bR^{3}$ is Alexandrov if it extends to an 
immersion of the bounding handlebody $M$). Similar remarks hold for all $\k \in \bR$.  

  Let $\cC^{m.\a}$ denote the space of conformal classes $[\g]$ of $C^{m,\a}$ metrics 
$\g$ on $\dm$, and let $H$ denote the mean curvature of $\dm \subset (M, g)$, with 
respect to the outward unit normal. It is proved in \cite{An2} that the moduli space 
$\cM_{\k} = \cM_{\k}^{m,\a}$ is a ($C^{\infty}$) smooth Banach manifold. Moreover, 
setting as above $\gamma = g|_{T(\dm)}$, the map
\begin{equation}\label{1.1}
\Pi: \cM_{\k} \to \cC^{m,\a}(\dm) \times C^{m-1,\a}(\dm),
\end{equation}
$$\Pi([g]) = ([\g], \ H),$$
is a ($C^{\infty}$) smooth Fredholm map, of Fredholm index 0. In fact the boundary data 
in \eqref{1.1} form elliptic boundary data for the Einstein equations. There are other 
elliptic boundary value problems for Einstein metrics, some of which are discussed 
in \cite{An2}. However, the data in \eqref{1.1} is geometrically the most natural so we 
restrict the discussion to this case. 

\medskip

  It follows in particular that $Im \Pi$ is a variety of finite codimension in 
$\cC^{m,\a}(\dm) \times C^{m-1,\a}(\dm)$. One would expect that generic metrics 
in $\cM_{\k}$ are regular points for $\Pi$, in which case $Im \Pi$ would at least 
contain open domains in the target $\cC^{m,\a}(\dm) \times C^{m-1,\a}(\dm)$; (a proof 
of this is still lacking however).

  The result in \eqref{1.1} shows that one has a good local existence theory for this 
boundary value problem and it raises the global problem:

\medskip

{\sc Question}. Given $([\g], H) \in \cC^{m,\a}(\dm) \times C^{m-1,\a}(\dm)$, (possibly with 
some restrictions), does there exist a unique metric, (up to isometry), $g$ on $M$ such 
that 
\be \label{1.2}
\Pi(g) = ([\g], H).
\ee

  To the author's knowledge, it does not seem that this question, although clearly quite 
natural, has been studied previously. There are many previous works on the existence of 
surfaces of prescribed mean curvature in $\bR^{3}$ for instance, cf.~\cite{TW} for example, 
and \cite{Y} or \cite{GT} for further discussion and references. However, in these 
situations $H$ is a given function on $\bR^{3}$; moreover, there is no prescription 
of the conformal class. In the case of $\k = 0$ for example, the question can be rephrased 
as the question of the existence and uniqueness of an Alexandrov immersion of a surface 
$F: \Sigma = \dm \to \bR^{3}$ with prescribed conformal class $[\g]$ and prescribed 
mean curvature $H$, i.e.
\be \label{1.3}
[F^{*}(g_{Eucl})] = [\g], \ \ H(F(x)) = H(x).
\ee
Note that the diffeomorphism group of $\dm$ acts non-trivially on both parts $[\g]$ 
and $H$ of the boundary data: if $\phi \in \Diff(\dm)$, then the immersion $F\circ 
\phi$ has the same image as $F$, but is a reparametrization of $F$. For different but 
related studies on surfaces of prescribed mean curvature, see for example \cite{Ka}, 
\cite{Ke} and \cite{LT}. 
 
\medskip

   To address a global question as above, the basic issue is whether the boundary map $\Pi$ 
is proper. In analogy to the simpler method of continuity commonly used in PDE, this is 
the closedness issue; one requires apriori estimates or compactness properties for spaces 
of solutions. For very simple reasons, the map $\Pi$ is not proper in general, and one 
first needs to sharpen the problem to account for this. Thus, for example smooth bounded 
domains in $\bR^{3}$ may degenerate from the ``inside'', in that the injectivity radius 
within $M$ may go to 0 near $\dm$, causing the boundary to develop self-intersections and 
the domain $M$ is no longer a manifold. This behavior can be ruled out via the maximum 
principle, (see also Lemma 2.4 below), under the assumption that $H > 0$. Thus, let
$$\cM_{\k}^{+} = \Pi^{-1}(\cC^{m,\a} \times C_{+}^{m-1,\a}),$$
be the space of constant curvature metrics with $H > 0$ at $\dm$. Clearly, $\cM_{\k}^{+}$ 
is an open submanifold of $\cM_{\k}$ and one may consider the associated (restricted) 
boundary map 
\begin{equation} \label{1.4}
\Pi_{+}: \cM_{\k}^{+} \to \cC^{m,\a} \times C_{+}^{m-1,\a}.
\end{equation}

  Next, recall by the uniformization theorem that the space of metrics $Met(S^{2})$ on 
$S^{2}$ equals $\Diff(S^{2})\times C_{+}(S^{2})$; any metric $\g$ is of the form $\g = 
\phi^{*}(\l^{2}\g_{+1})$, where $\g_{+1}$ is the round metric of radius $1$. (Similarly 
for surfaces $\Sigma$ of higher genus, $Met(\Sigma)$ is a bundle over the Riemann 
moduli space with fiber $\Diff(\Sigma)\times C_{+}(\Sigma)$). The group $\Diff(S^{2})
\times C_{+}(S^{2})$ thus acts transitively on $Met(S^{2})$ and has stabilizer at 
$\g_{+1}$ equal to the group of essential conformal transformations $\Conf(S^{2})$ 
of $S^{2}(1)$. It follows that one has a natural identification 
\be \label{1.5}
\cC^{m,\a}(S^{2}) \simeq \Diff^{m+1,\a}(S^{2})/\Conf(S^{2}).
\ee
The conformal group also acts on the space $C_{+}^{m-1,\a}$ of mean curvature 
functions: $H \to H\circ \f$, for $\f \in \Conf(S^{2})$. It is easy to verify that 
this action is free and proper, except on the functions $H = const$, which are 
the fixed points of the action; this is because the flow of the conformal vector 
fields contracts or expands all of $S^{2}\setminus \{pt\}$ to a point. 

  It follows then that at the special values $([\g], c)$ where $H = c$, the 
map $\Pi_{+}$ in \eqref{1.4} is not proper. The non-compact conformal group 
$\Conf(S^{2})$ fixes this data, but acts nontrivially (and faithfully) on 
$\cM_{\k}^{+}$; if $\Pi_{+}(g) =  ([\g], c)$, then also $\Pi_{+}(\f^{*}(g)) 
= ([\g], c)$, for any $\f \in \Conf(S^{2})$ extended to a diffeomorphism of $M$. 
On the other hand, this is the only value where $\Conf(S^{2})$ acts non-properly. 
(Note this issue arises only for $S^{2}$, not for boundaries of higher genus). 

   There are two ways to deal with this issue. First, one may just study the 
behavior of $\Pi_{+}$ away from the ``round'' metrics $H = c$, i.e.~consider 
the global behavior of the map
\begin{equation} \label{1.6}
\Pi': \cM_{\k}^{'} \to \cC^{m,\a} \times (C_{+}^{m-1,\a})',
\end{equation}
where $(C_{+}^{m-1,\a})' = C_{+}^{m-1,\a} \setminus \{constants\}$ and 
$\cM_{\k}^{'} = \Pi_{+}^{-1}((C_{+}^{m-1,\a})')$. This map is again smooth 
and Fredholm, of index 0. 

  Alternately, one may include the round metrics, but divide out by the action 
of $\Conf(S^{2})$. Briefly, as is standard, choose a fixed marking to freeze the 
action of the conformal group on $S^{2}$. Thus, fix three points $p_{i}$, 
$i = 1,2,3$ on $S^{2}(1)$ with 
\be \label{1.7}
dist_{\g_{+1}}(p_{i}, p_{j}) = \pi/2.
\ee
Let $\cN_{\k}^{+}$ be the marked moduli space of constant curvature metrics on $M = 
B^{3}$ consisting of metrics $g$ such that \eqref{1.7} holds on $\dm = S^{2}$ with 
$\g = g|_{\dm}$ in place of $\g_{+1}$. The condition \eqref{1.7} can always be 
realized by changing $g$ by a conformal diffeomorphism, so that 
$$\cN_{\k}^{+} \simeq \cM_{\k}^{+}/\Conf(S^{2}),$$
the condition \eqref{1.7} giving a slice for the action of $\Conf(S^{2})$ on 
$\cM_{\k}$. (Of course this marking is not necessary in case $\chi(\dm) \leq 0$ 
so that $\cN_{\k}^{+} \simeq \cM_{\k}^{+}$ in such cases). The map $\Pi_{+}$ in 
\eqref{1.4} clearly restricts under the slicing \eqref{1.7} to a smooth map 
\begin{equation}\label{1.8}
\Pi_{+}: \cN_{\k}^{+} \to \cC^{m,\a}(\dm) \times C^{m-1,\a}(\dm).
\end{equation}
However, the index of this map is now -3, since $\Conf(S^{2})$ is 3-dimensional and 
so one must also divide the target space by the remaining action of $\Conf(S^{2})$ on 
the space $C_{+}^{m-1,\a}$ of mean curvature functions. Let then 
$\cB[\cC^{m,\a} \times C_{+}^{m-1,\a}]$ be the quotient space 
$$\cB[\cC^{m,\a} \times C_{+}^{m-1,\a}] = [\cC^{m,\a} \times C_{+}^{m-1,\a}]/\Conf(S^{2}) 
= \cC^{m,\a} \times (C_{+}^{m-1,\a}/\Conf(S^{2})).$$
The map $\Pi_{+}$ in \eqref{1.8} then descends to a smooth Fredholm map
\begin{equation}\label{1.9}
\Pi_{+}: \cN_{\k}^{+} \to \cB[\cC^{m,\a} \times C_{+}^{m-1,\a}],
\end{equation}
of Fredholm index 0. The two formulations \eqref{1.6} and \eqref{1.9} are basically 
equivalent. 

  Now if $\Pi_{+}$ in \eqref{1.9}, (or $\Pi'$ in \eqref{1.6}), is proper, then by work 
of Smale \cite{Sm}, it has a well-defined degree $deg \Pi' \in \bZ_{2}$, (and most 
likely a $\bZ$-valued degree if the spaces can be given an orientation). Elementary 
degree theory implies that if 
$$deg \Pi_{+} \neq 0,$$
then $\Pi_{+}$ in \eqref{1.6} is surjective, answering at least the existence part of 
the question above. 

  In fact, if $\Pi_{+}$ is proper, then one has
\be \label{1.11}
deg \Pi_{+} \neq 0 \ \ {\rm for} \ \ \dm = S^{2}, \ \ {\rm but} \ \  deg \Pi_{+} = 0 \ \ 
{\rm for} \ \ \dm \neq S^{2}.
\ee
This follows from the Alexandrov-Hopf rigidity theorems, \cite{Al}, \cite{H}. Namely, 
any metric on a surface $\Sigma = \Sigma_{g}$ of genus $g$ which is Alexandrov immersed 
in a space-form with $H = const$ is necessarily a round sphere. This uniqueness also holds 
infinitesimally, showing that the ``round'' conformal class $([\g_{+1}], H = c)$ 
is a regular value of $\Pi$. Since the Hopf theorem implies that the inverse 
image of this round regular value is unique, it follows that $deg \Pi_{+} \neq 0$ 
for $\Sigma = S^{2}$. For $\Sigma_{g}$ with $g \neq 0$, the same argument shows 
that $\Pi_{+}$ in \eqref{1.9} is not surjective, which implies $deg \Pi_{+} = 0$. 

  Thus, at least in the case of $S^{2}$ the existence question above has been 
reduced to the properness of $\Pi_{+}$. This issue will be discussed in detail 
in the next section; we will show however that $\Pi_{+}$ or $\Pi'$ is in fact 
not proper, so that further modifications are necessary to understand the global 
behavior of these boundary maps.

\section{Analysis of the Boundary map $\Pi$.}

  We begin by filling in some details from the discussion in \S 1. Since the 
full curvature is determined by the Ricci curvature in 3-dimensions, any metric 
$g \in \bM_{\k}$ satisfies the Einstein equation
\be \label{2.1}
Ric_{g} - 2\k\cdot g = 0.
\ee
We wish to view \eqref{2.1} as an elliptic equation for $g$. This is not possible due 
to the diffeomorphism invariance of \eqref{2.1}, and so one needs to choose a gauge to 
break this invariance. Let $\widetilde g \in \bM_{\k}$ be a fixed but arbitrary 
(constant curvature) background metric. The simplest choice of gauge is the 
Bianchi-gauge, with the associated Bianchi-gauged Einstein operator, 
given by 
\begin{equation} \label{2.2}
\Phi_{\widetilde g}: Met(M) \rightarrow  S_{2}(M), 
\end{equation}
$$\Phi_{\widetilde g}(g) = Ric_{g} - 2\k g + 
\delta_{g}^{*}\beta_{\widetilde g}(g),$$
where $(\delta^{*}X)(A,B) = \frac{1}{2}(\langle \nabla_{A}X, B\rangle  + 
\langle \nabla_{B}X, A\rangle )$ and $\delta X = -tr \delta^{*}X$ is the 
divergence and $\beta_{\widetilde g}(g) = \delta_{\widetilde g} g + 
\frac{1}{2}d tr_{\widetilde g} g$ is the Bianchi operator with respect to 
$\widetilde g$. 

 Clearly $g$ is Einstein if $\Phi_{\widetilde g}(g) = 0$ and 
$\beta_{\widetilde g}(g) = 0$, so that $g$ is in the Bianchi-free gauge with 
respect to $\widetilde g$. Using standard formulas for the linearization 
of the Ricci and scalar curvatures, cf.~\cite{B} for instance, one finds 
that the linearization of $\Phi$ at $\widetilde g = g$ is given by 
\begin{equation} \label{2.3}
L(h) = 2(D\Phi_{\widetilde g})_{g}(h) = D^{*}Dh - 2R(h). 
\end{equation}
The zero-set of $\Phi_{\widetilde g}$ near $\widetilde g$, 
\begin{equation} \label{2.4}
Z = \{g: \Phi_{\widetilde g} = 0\}, 
\end{equation}
consists of metrics $g \in Met(M)$ satisfying the equation $Ric_{g} - 2\k g + 
\delta_{g}^{*}\beta_{\widetilde g}(g) = 0$. 

  Given $\widetilde g$, consider the Banach space
\begin{equation}\label{2.5}
Met_{C}(M) = Met_{C}^{m,\alpha}(M) = \{g \in Met^{m,\alpha}(M): 
\beta_{\widetilde g}(g) = 0 \ {\rm on} \ \partial M \}. 
\end{equation}
Clearly the map 
$$\Phi : Met_{C}(M) \rightarrow  S^{2}(M),$$
is $C^{\infty}$ smooth. Let $Z_{C}$ be the space of metrics $g\in Met_{C}(M)$ 
satisfying $\Phi_{\widetilde g}(g) = 0$, and let $\bM_{C} = \bM_{\k} \cap Z_{C}$
be the subset of constant curvature metrics $g$, $Ric_{g} = 2\k g$ in 
$Z_{C}$. It is proved in \cite{An2} that $Z_{C}$ is a smooth Banach manifold and 
$$\bM_{C} = Z_{C},$$
so that any metric $g\in Z_{C}$ near $\widetilde g$ is necessarily constant 
curvature, with $Ric_{g} = 2\k g$, and in Bianchi gauge with respect to 
$\widetilde g$. This result also holds at the linearized level. The spaces 
$Z_{C}$ are smooth slices for the action of the diffeomorphism group 
$\Diff_{1}^{m+1,\a}$ on $\bM_{\k}$ and it follows that the quotient 
$\cM_{\k}$ is a smooth Banach manifold.

\medskip

  Next consider elliptic boundary data for the operator $\Phi$ in \eqref{2.2}. Dirichlet 
or Neumann boundary data are not elliptic; this follows by inspection from the Gauss 
constraint equation \eqref{2.12} below, (or from the proof below). The following result 
is proved in \cite{An2}; we give the main details of the proof, since it is useful to 
compare this with the discussion in \S 3.

\begin{proposition}\label{p2.1}
The Bianchi-gauged Einstein operator $\Phi$ with boundary conditions
\begin{equation}\label{2.6}
\beta_{\widetilde g}(g) = 0, \ \ [g^{T}] = [\gamma], \ \ 
H_{g} = h \ \ {\rm at} \ \ \partial M,
\end{equation}
is an elliptic boundary value problem of Fredholm index 0.
\end{proposition} 

{\sc Proof:} It suffices to show that the leading order part of the linearized operators 
at the Euclidean metric forms an elliptic system. The leading order symbol of $L = D\Phi$ 
is given by 
\begin{equation}\label{2.7}
\sigma(L) = -|\xi|^{2}I,
\end{equation}
where $I$ is the $3\times 3$ identity matrix. In the following, 
the subscript 0 represents the direction normal to $\partial M$ in $M$, and Latin 
indices run from $1$ to $2$. The positive roots of \eqref{2.7} are $i|\xi|$, with 
multiplicity $3$. Writing $\xi = (z, \xi_{i})$, the symbols of the leading order terms 
in the boundary operators are:
$$-2izh_{0k} - 2i\sum \xi_{j}h_{jk} + i\xi_{k}tr h = 0,$$
$$-2izh_{00} - 2i\sum \xi_{k}h_{0k} + iztr h = 0,$$
$$h^{T} = (\gamma')^{T} \ \ mod\, \gamma, \ \  H_{h}' = \omega,$$
where $h$ is a $3\times 3$ matrix. Ellipticity requires that the operator defined 
by the boundary symbols above has trivial kernel when $z$ is set to the root $i|\xi|$. 
Carrying this out then gives the system 
\begin{equation}\label{2.8}
2|\xi|h_{0k} - 2i\sum \xi_{j}h_{jk} + i\xi_{k}tr h = 0,
\end{equation}
\begin{equation}\label{2.9}
2|\xi|h_{00} - 2i\sum \xi_{k}h_{0k} - |\xi|tr h = 0,
\end{equation}
\begin{equation}\label{2.10}
h_{kl} = \phi\delta_{kl}, \ \ H_{h}' = 0.
\end{equation}
where $\phi$ is an undetermined function. 

  Multiplying \eqref{2.8} by $i\xi_{k}$ and summing gives
$$2|\xi|i\sum \xi_{k}h_{0k} =  2i^{2}\xi_{k}^{2}h_{kk} - i^{2}\xi_{k}^{2}tr h.$$
Substituting \eqref{2.9} on the term on the left above then gives
$$2|\xi|^{2}h_{00} - |\xi|^{2}tr h =  -2\sum \xi_{k}^{2}h_{kk} + |\xi|^{2}tr h,$$
so that
$$|\xi|^{2}h_{00} - |\xi|^{2}tr h =  -\sum \xi_{k}^{2}h_{kk} = 
-\phi |\xi|^{2}.$$
Using the fact that $tr h - h_{00} = \sum h_{kk} = n\phi$, it follows that 
$\phi = 0$ and hence $h^{T} = 0$. 

  A simple computation shows that to leading order, $H_{h}' = tr^{T}(\nabla_{N}h - 
2\delta^{*}(h(N)^{T}))$, which has symbol $iz\sum h_{kk} - 2i\xi_{k}h_{0k}$. Setting 
this to 0 at the root $z = i|\xi|$ gives
$$\sum (|\xi|h_{kk} + 2i\xi_{k}h_{0k}) = 0.$$
Since $h^{T} = 0$, this gives $\sum \xi_{k}h_{0k} = 0$, which, via \eqref{2.9} 
gives $h_{00} = 0$ and hence via \eqref{2.8}, $h = 0$. 

  This proves that the boundary data \eqref{2.6} are elliptic for $\Phi$. The proof that 
the Fredholm index is 0 is given in \cite{An2}. 

{\endproof}

  We now turn to the main issue, the properness of the map $\Pi_{+}$ in \eqref{1.9} 
or $\Pi'$ in \eqref{1.6}. This amounts to proving (apriori) estimates for metrics 
$g \in \cM_{\k}$ in terms of the boundary data $([\g], H)$. The main result in 
this direction is the following:
\begin{proposition}\label{p2.2}
Let $\cK$ be a compact set in the space of boundary data $\cC^{m,\a} 
\times C_{+}^{m-1,a}$. Then for any $K < \infty$, the space of metrics 
$g \in \cN_{\k}^{+}$ such that 
\be \label{2.11}
\Pi_{+}(g) \in \cK\ \ {\rm and} \ \ a = area(\dm) \leq K,
\ee
is compact. 
\end{proposition}

  This result shows that $\Pi_{+}$ is proper, under the assumption of an upper bound 
on $a$. The proof of this result follows below, organized into several lemmas. 

\medskip

  To begin, we recall the constraint equations at $\dm$, i.e.~the Gauss and 
Gauss-Codazzi equations:
\be \label{2.12}
|A|^{2} - H^{2} + 2K_{\g} = s_{g} - 2Ric_{g}(N,N) = 2\k,
\ee
\be \label{2.13}
\d(A - H\g) = -Ric(N, \cdot) = 0,
\ee
where $A$ is the second fundamental form and $N$ is the outward unit normal. 

  One of the most important issues is to obtain a bound on $|A|$. 
\begin{lemma}\label{l2.3}
There is a constant $C_{0} < \infty$, depending only on $\cK$ and $K$ in 
Proposition 2.2, such that
\be \label{2.14}
|A| \leq C_{0}.
\ee
\end{lemma}

{\sc Proof:} The proof is by contradiction, by means of a blow-up argument. To begin, 
integrating the Gauss constraint \eqref{2.12} and using the Gauss-Bonnet theorem gives
\be \label{2.15}
\int_{\dm}|A|^{2} = \int_{\dm}H^{2} - 4\pi \chi(\dm) + 2\k \cdot area(\dm).
\ee
By assumption, $a = area(\dm)$ is uniformly bounded, 
$$area(\dm) \leq K < \infty.$$
This and \eqref{2.15} give an apriori bound on the scale-invariant quantity $\int |A|^{2}$: 
\be \label{2.16}
\int_{\dm} |A|^{2} \leq C.
\ee

   Now choose a point $x$ on $\dm$ where $|A|$ is maximal, and rescale the metric so that 
$|A|(x) = 1$, with $|A|(y) \leq 1$ everywhere, so $\bar g = \l^{2}g$ where $\l = |A|(x)$. 
It follows directly from the constraint equation \eqref{2.12} that the intrinsic curvature 
$K_{\bar \g}$ of $\bar \g$ is also uniformly bounded. The family of such metrics is 
compact in the pointed $C^{1,\a}$ topology, by the Cheeger-Gromov compactness theorem for 
instance; this means that modulo diffeomorphisms of $\dm$, the metric $\g$ itself 
is uniformly controlled in $C^{1,\a}$, (in suitable local coordinates and within 
bounded distance to $x$). 

   Now by assumption, the conformal class $[\g]$ of $\g$ is uniformly controlled. It 
follows that the diffeomorphisms above, (in which $\bar \g$ is uniformly controlled), 
are themselves controlled modulo the group of conformal diffeomorphisms, cf.~\eqref{1.5}. 
Thus, passing if necessary from $\bar \g$ to $\g' = \phi^{*}(\bar \g)$, where $\phi$ is a 
conformal diffeomorphism, it follows that the metric $\g'$ is uniformly controlled 
in $C^{1,\a}$, (locally, within bounded distance to $x$). Together with the uniform 
bound on $|A|$ above, it follows from Proposition 2.1 and elliptic regularity that 
the metric $g'$ is controlled in the stronger $C^{m,\a}$ norm, up to its boundary. 

  Suppose then $g_{i}$ is a sequence where $\max |A| \to \infty$. By rescaling as above 
one may pass to a smoothly convergent subsequence of $\{g_{i}'\}$ to obtain a smooth 
limit $g'$. The smooth ($C^{m,\a}$) convergence implies on the one hand that the 
limit is not flat, since the condition $|A|(x) = 1$ passes continuously to the limit. 
The estimate \eqref{2.16} also holds on the limit. Since $H \to 0$ in the rescalings, 
it follows that the limit is a complete immersed minimal surface in $\bR^{3}$ with 
finite total curvature 
$$\int_{\Sigma}|A|^{2} < \infty.$$
Moreover, since the conformal classes $[\g_{i}]$ of $g_{i}$ on $\dm$ are uniformly 
controlled, the sequence $\g_{i}'$ has a uniformly controlled (large scale) atlas of 
conformal coordinates. Hence the limit is conformally isometric to $\bR^{2}$, 
i.e.~the limit minimal surface is pointwise conformal to $\bR^{2}$. Finally, these 
minimal surfaces are in fact embedded; this follows from Lemma 2.4 below. However, 
it is well-known that the only such surfaces are flat planes, (cf.~\cite{MRR} for 
instance), and hence $A = 0$ in the limit. This contradiction establishes the 
bound \eqref{2.14}. 

{\endproof} 

  Next we show that the normal exponential map has injectivity radius bounded below
$$inj_{N} \geq i_{0}$$
where $i_{0}$ depends only on an upper bound for $|A|$. This follows from the following 
Lemma. 

  Let here $N$ be the inward unit normal to $\dm$ in $M$ and consider the associated normal 
exponential map to $\dm$, $tN \to \exp_{p}(tN)$, giving the geodesic normal to $\dm$ at 
$p$. This is defined for $t$ small, and let $D(p)$ be the maximal time interval on which 
$\exp_{p}(tN) \in M$, (so that the geodesic does not hit $\dm$ again before time $D(p)$). 
Thus, $D: \dm \to \bR^{+}$. 

\begin{lemma}\label{l2.4}
Given $H > 0$, suppose $|A| \leq C_{0}$. Then there is a constant $t_{0}$, depending only on 
$C_{0}$ and $\k$, (and the lower bound on $H$ when $\k < 0$), such that 
\be \label{2.17}
D(p) \geq t_{0}.
\ee
\end{lemma}

{\sc Proof:} This is a well-known result in Riemannian geometry, essentially due to Frankel, 
and follows from the 2nd variational formula for geodesics. First, given bounds on $|A|$ and 
$\k$, by standard comparison geometry one has a lower bound on the distance to the focal 
locus of the normal exponential map $\exp(tN)$, i.e.~a lower bound $d_{0}$ on the distance 
to focal points. Suppose then 
$$\min D < d_{0}.$$
If the minimum is achieved at $p$, then the normal geodesic to $\dm$ at $p$ intersects 
$\dm$ again at a point $p'$, and the intersection is orthogonal to $\dm$ at $p'$. Denoting 
this geodesic by $\s$, and letting $\ell = D(p)$ be the length of $\s$, the 2nd variational 
formula of energy gives
\be \label{2.18}
E''(V,V) = \int_{0}^{\ell}(|\nabla_{T}V|^{2} - \langle R(T,V)V,T\rangle) dt - 
\langle \nabla_{V}T, V \rangle|_{0}^{\ell},
\ee
where $T = \dot \s$ and $V$ is any variation vector field along $\s$ orthogonal to $\s$. 
By the minimizing property, one has $E''(V,V) \geq 0$, for all $V$. Choose then $V = 
V_{i}$ to be parallel vector fields $e_{i}$, running over an orthonormal basis at 
$T_{p}(\dm)$. The first term in \eqref{2.18} then vanishes, while the second sums to 
$-Ric(T,T) = -2\k$. The boundary terms sum to $\pm H$, at $p$ and $p'$. Taking into 
account that $T$ points into $M$ at $p$ while it points out of $M$ at $p'$, this gives
$$0 \leq -2\k\ell - (H(p) + H(p')).$$
Since $H > 0$, this gives immediately a contradiction if $\k \geq 0$, and also gives a 
contradiction if $\k < 0$, if $H$ is bounded below, depending only on the size of $\k$ 
(if $\ell$ is sufficiently small). This proves the estimate \eqref{2.17}. 

{\endproof}

  Finally we show that Lemma 2.3 implies that the intrinsic geometry of $(\dm, \g)$ is 
controlled. 
\begin{lemma}\label{l2.5}
There is a constant $C_{1}$, depending only on $C_{0}$ in \eqref{2.14} and $\cK$, 
such that
\be \label{2.20}
|K_{\g}| \leq C_{1},
\ee
Moreover, the metric $\gamma$ is uniformly controlled, modulo conformal 
diffeomorphisms, by $C_{1}$, (and $a$). 
\end{lemma}

{\sc Proof:}  As in the proof of Lemma 2.3, via the Gauss constraint equation 
\eqref{2.12}, a bound on $|A|$ implies a bound on $K_{\g}$, giving \eqref{2.20}. 
A standard simple analytic argument then gives control on the metric $\g$ itself 
when $\chi(\dm) \leq 0$. Namely, write $\g = \lambda^{-2}\g_{0}$, where 
$\g_{0}$ is the conformal metric with constant curvature $\sigma$ and $\sigma$ is 
chosen so that $area \g = area \g_{0}$. The formula for the behavior of Gauss 
curvature under conformal changes then gives
$$\l^{2}\D_{\g_{0}}(\log \l) = -\l^{2}\sigma - K_{\g}.$$
The maximum principle implies an upper bound on $\l$ and hence, by elliptic regularity, 
one has uniform $C^{1,\a}$ control on $\l$ and so, (via standard bootstrap arguments), 
$\l$ is controlled in $C^{m,\a}$. 

  This argument does not work when $\chi(\dm) > 0$, (since the minimum or maximum 
principle does not hold). In this case, one can use the same argument as that given 
in the proof of Proposition 2.2. Thus, the bound $|K_{\g}|$ implies that the 
metric is controlled modulo diffeomorphisms, by the Cheeger-Gromov compactness 
theorem. Here we use the fact that the length of the shortest closed geodesic, 
and hence the injectivity radius of $\g$, is bounded below, since $|A|$ is bounded 
above. Since the conformal class $[\g]$ is assumed to be controlled, the 
diffeomorphisms are themselves controlled, modulo the group of conformal 
diffeomorphisms. 

  Note that \eqref{2.15} shows that $a = area(\dm)$ is bounded below, and hence 
the diameter of $(\dm, \g)$ is also bounded above and below. This proves the 
result. Note also that since $diam (M, g) \leq diam(\dm, \g)$, this also 
gives a uniform upper bound on the diameter of $(M, g)$. 

{\endproof}

  The results above prove Proposition 2.2. This result implies that the ``enhanced'' 
boundary map
\be \label{2.21}
g \to ([\g], H, a)
\ee
is proper. While this map is Fredholm, it is Fredholm of index -1 and so does not have 
a well-defined degree; for this, one needs the Fredholm index to be non-negative. 

  There are two ways in which one may try to proceed at this point. 

\medskip

{\bf (I).} One may try to prove that $a = area(\dm)$ is controlled by the 
boundary data $([\g], H)$, which would then prove that $\Pi$ itself, 
(i.e.~$\Pi_{+}$ or $\Pi'$), is proper. 

  However, this is false. It follows from the proof of Proposition 2.2 that 
counterexamples must closely resemble the helicoid (in a suitable scale), since 
the helicoid is the unique complete embedded minimal surface in $\bR^{3}$ 
conformally equivalent to $\bR^{2}$, (besides the plane), cf.~\cite{MR}. In 
fact conversely, one may use the helicoid to construct examples of metrics 
$g_{i}$ where $([\g_{i}], H_{i})$ are uniformly bounded but 
\be \label{2.21a}
a_{i} \to \infty.
\ee
To see this, consider the helicoid ${\mathcal H} = {\mathcal H}_{L}$, 
$$x = \rho \cos\theta, \ \ y = \rho \sin \theta, \ \ z = L^{-1}\theta,$$
$L = L(\e) >> 1$, wrapping around $z$-axis arbitrarily many times in the interval 
$z \in [-\e, \e]$; assume here $\rho \in [-1,1]$. Consider an almost horizontal 
$S^{1} \subset {\mathcal H}$ formed by connecting two line segments parallel to 
the $x$-axis in ${\mathcal H}$ at height $\pm \e$ by a circular arc joining their 
endpoints along a helix in ${\mathcal H}$. This $S^{1}$ bounds a disc $D^{2} \subset 
{\mathcal H}$. Now form the vertical $z$-cylinder over this boundary $S^{1}$, and 
take it to a fixed height, say $z = \frac{1}{2}$ and then cap off the circular 
boundary at $z = \frac{1}{2}$ by a horizontal disc. This gives first an immersed 
$S^{2}$, which is also Alexandrov immersed, since it may be perturbed to an 
embedding. This $S^{2}$ may also be perturbed so that $H > 0$ everywhere. Namely, 
one may first deform the helicoid very slightly to a surface with $H > 0$, and 
in fact with $H$ uniformly bounded away from 0 and $\infty$, cf.~\cite{Ti} for 
instance. The vertical cylinder has $H > 0$ and one can bend the top flat disc 
to $H > 0$. Finally, the corners of $S^{2}$ may also be smoothed to $H > 0$ 
everywhere. 

  The conformal class of the helicoid is fixed under arbitrary rescalings, 
(i.e.~variations of $L$), and the gluing process above is also uniformly 
controlled; hence the conformal class of the collection of surfaces above is 
uniformly controlled. It is clear from the construction that \eqref{2.21a} holds 
as the number of wrappings of the helicoid is taken to infinity. 

  Recall the Hopf uniqueness theorem: if $\Sigma$ is a sphere immersed in a 
space-form of constant curvature with $H = const$, then $\Sigma$ is umbilic and 
so locally isometric to a round sphere. The examples above seem to indicate or 
suggest that the rigidity associated to the Hopf theorem cannot be weakened to 
an ``almost rigidity'' theorem; thus we expect given any $\e > 0$, there exist 
surfaces $\Sigma_{\e} \subset \bR^{3}$ diffeomorphic to $S^{2}$ such that 
$$2 - \e \leq H_{\Sigma_{\e}} \leq 2 + \e,$$ 
which are not close to a round sphere $S^{2}(1) \subset \bR^{3}$. Of course, 
one must have $area (\Sigma_{\e}) \to \infty$ as $\e \to 0$. It would be interesting 
to know the answer to this question.

\medskip

{\bf (II).} Instead, refering to the context of \eqref{2.21}, one may add an extra 
scalar variable $\l$ to the domain to obtain a map of Fredholm index 0. [Alternately, 
one may restrict the domain $\cM_{\k}'$ in \eqref{1.6}, or $\cN_{\k}^{+}$ in \eqref{1.9} 
to surfaces where $a = 1$; correspondingly, one must then divide the target space by an 
$\bR^{+}$ action. There is no essential difference between these so we discuss only 
the former]. 

  Thus, extend for example the domain $\cM_{\k}'$ to $\cM_{\k}'\times \bR^{+}$, and 
consider the following typical examples:
\be \label{2.22}
(g, \l) \to ([\g], H, a + \l).
\ee
\be \label{2.23}
(g, \l) \to ([\g], \frac{H}{H_{\min}} - 1 + \l, a),
\ee
where $H_{\min}$ is the minimum value of $H$. 

  These maps are Fredholm, of Fredholm index 0. However, neither map is proper; 
in \eqref{2.22}, one may have $\l \to 0$ while in \eqref{2.23} one may have 
$H_{\max} \to \infty$, both within compact sets of boundary data. Consider next 
shifting the $\l$-variable also to the space $\cC$. For example, let 
$\psi_{\l}$ be a curve of diffeomorphisms of $\dm$ with $\psi_{\l} \to \infty$ 
as $\l \to 0$ and consider
\be \label{2.24}
(g, \l) \to ([\psi_{\l}^{*}(\g)], H, a + \l).
\ee
As before, this behaves well in the second two factors, but it is not clear if there 
exists a curve $\psi_{\l}$ for which \eqref{2.24} is proper; in any case we have not 
been able to find a construction to make this map proper. 

  Next consider enlarging the domain by adding a scale factor, when $\k \neq 0$. Thus, 
let $\cM' = \cup_{\k < 0}\cM_{\k}'$, and define
\be \label{2.25}
\cM' \to \cC\times C_{+},
\ee
$$g \to ([\g], H, a).$$
This gives of course control of both $H$ and $a$, and by \eqref{2.15}, one deduces 
uniform control on $|\k|$, when $\k < 0$. However, one cannot prevent the possibilty 
that $\k \to 0$, so that again its not clear if this map can be made proper. 
(Including the spaces $\cM_{\k}'$ with $\k > 0$ also does not seem to help). 

  Consider finally the following modification of \eqref{2.23}:
\be \label{2.26}
\widetilde \Pi: (g, \l) \to ([\g], \frac{H}{H_{\min}} - 1 + \l, a + H_{\max}).
\ee
The map $\widetilde \Pi$ is Fredholm, of index 0, and is now proper by Proposition 2.2, 
since control of the data in the target space gives control on $H$, $a$ and the 
conformal class $[\g]$. It thus has a well-defined degree. 

  However, the Hopf rigidity theorem now shows that 
\be \label{2.27}
deg \widetilde \Pi = 0.
\ee
Namely, consider the case $\dm = S^{2}$ and $\k = 0$. When restricted to the 
``round'' metrics where $H = const$, by the constraint equation \eqref{2.12} one has
$$H^{2}a = 16\pi,$$
so that $H = H_{\max} = \sqrt{\frac{16\pi}{a}}$. This gives
$$\beta(a) \equiv a + H_{\max} = a + \sqrt{\frac{16\pi}{a}},$$
which, for a given value of $\beta(a) = c$ has two positive real solutions 
$a > 0$. The function $\beta$ is a simple fold map $\bR^{+} \to [(4\pi)^{1/3},\infty)$. 
This implies \eqref{2.27}. (Although the round metric is not in $\cM_{\k}'$, the 
discussion above remains valid for data near the round metric). 

  There is another, quite different argument showing that $\widetilde \Pi$ is not 
onto, and hence has degree 0. Namely on any $g \in \cM_{\k}$ with boundary data 
$(\g, H)$, one has
\be \label{2.28}
\int_{S^{2}}X(H)dV_{\g} = 0,
\ee
where $X$ is any conformal Killing field on $S^{2}$. This follows from the constraint 
equation \eqref{2.13}. Namely, pairing \eqref{2.13} with a vector field $X$ and 
integrating over $(\dm, \g)$ gives 
$$\int_{\dm}\langle \d A, X \rangle = -\int_{\dm}\langle dH, X \rangle.$$
The left side equals $\int_{\dm}\langle A, \d^{*}X\rangle$, and for $X$ conformal 
Killing, $\int_{\dm}\langle A, \d^{*}X\rangle = \frac{1}{n}\int_{\dm}H div X$. 
On the other hand, $-\int_{\dm}\langle dH, X \rangle = \int_{\dm}H div X$, which gives 
\eqref{2.28}, since $dim \dm = n \neq 1$. The result \eqref{2.28} is essentially due to 
\cite{Am}, although the proof given here is much simpler. The condition \eqref{2.28} is 
of course reminiscent of the Kazdan-Warner type obstruction \cite{KW} for the prescribed 
Gauss curvature problem. As in the Gauss curvature problem, note that the condition 
\eqref{2.28} is not conformally invariant, i.e.~it is not a well-defined condition 
on the target space $\cC^{m,\a} \times C_{+}^{m-1,\a}$. 

  The ``balancing condition'' \eqref{2.28} implies for instance that any $H$ which 
is a monotone function of a height function on $S^{2}(1)$ is not the mean curvature 
of a conformal immersion in a space-form. On the other hand, although \eqref{2.28} 
formally represents 3 independent conditions on $H$, it does not imply that 
$Im \Pi$ has codimension 3, (or any other codimension), in $\cC^{m,\a}\times 
C_{+}^{m-1,\a}$, again since it is not defined on this target space. 

\medskip

   Although we have succeeded in constructing a proper map $\widetilde \Pi$, it is 
not at all clear what $Im \widetilde \Pi$ is, or what the images of the closely related 
maps $\Pi_{+}$ and $\Pi'$ in \eqref{1.9}, \eqref{1.6} are. For instance, can $Im \Pi_{+}$ 
be described as the locus where a finite number of real-valued functions on the target 
are positive? Can one explicitly identify such functions characterising the boundary 
values of metrics in $\cM_{\k}$? 

  Finally, regardless of the surjectivity issue, the discussion in (I) above 
suggests that $\Pi_{+}$ is infinite-to-one, so highly non-unique. In sum, it 
would be interesting to understand these issues better, which seem on the 
whole much easier than the existence and uniqueness question for the 
isometric immersion problem discussed in \S 1.

\section{Generalization}

  Let $(M_{\k}, g_{\k})$ be any complete Riemannian 3-manifold of constant curvature $\k$; 
thus, up to scaling, $M_{\k}$ is one of $\bR^{3}$, $\bH^{3}$ or $\bS^{3}$ or a quotient 
of one of these spaces. Let $f: \dm \to M_{\k}$ be an Alexandrov immersion, and let $F$ 
denote an extension of $f$ to $M$. Then since $\pi_{1}(M, \dm) = 0$, the metric 
$F^{*}(g_{\k})$ is uniquely determined by the immersion $f$ on $\dm$, modulo $\Di_{1}(M)$. 
Thus, the map $\Pi_{+}$ in \eqref{1.4} is equivalent to a map
$$\Pi_{+}: Imm_{A}(\dm) \to  \cC^{m,\a}(\dm) \times C_{+}^{m-1,\a}(\dm).$$
This suggests that one could replace the space of metrics $\cM_{\k}$ by the space 
of immersions of $\Sigma = \dm$ into a space-form $M_{\k}$ or more generally 
into an arbitrary complete Riemannian manifold $(N, g_{N})$. This is in fact the case:
\begin{proposition}
Let $\Sigma = \Sigma_{g}$ be a compact surface of genus $g$ and let $(N, g_{N})$ be 
any complete Riemannian 3-manifold. Let $Imm^{m+1,\a}(\Sigma, N)$ be the space 
of $C^{m+1,\a}$ immersions of $\Sigma \to N$. Then the map
\be \label{3.2}
\Pi: Imm^{m+1,\a}(\Sigma, N) \to \cC^{m,\a}(\Sigma)\times C^{m-1,\a}(\Sigma),
\ee
$$\Pi(f) = ([f^{*}(g_{N})|_{\dm}], H(f(x))),$$
is a smooth Fredholm map of Fredholm index 0. 
\end{proposition}

{\sc Proof}: The space $Imm(\Sigma, N)$, is a smooth Banach manifold; the tangent 
space is given by the space of vector fields $v$ along a given immersion $f: \Sigma 
\to N$. The differential $D\Pi$ of $\Pi$ in \eqref{3.2} is given by
\be \label{3.3}
([(\d^{*}v)^{T}]_{0}, H_{\d^{*}v}'),
\ee
where $(\d^{*}v)^{T}$ is the restriction of $\d^{*}v$ to $T(\dm)$. The Fredholm property 
then follows by showing that the data \eqref{3.3} form an elliptic system of equations for 
$v$. We do this following the proof of Proposition 2.2. 

  Thus, write $v = v^{T} + fN$, where $v^{T}$ is tangent and $N$ is normal to 
$T(\dm)$. Then 
\be \label{3.4}
\d^{*}v = \d^{*}v^{T} + fA + df\cdot N,
\ee
so that $(\d^{*}v)^{T} = \d^{*}v^{T} + fA$. The second term is lower order in $v$ 
and so does not contribute to principal symbol. The principal symbol $\sigma$ of 
the first term is thus
\be \label{3.5}
\sigma([(\d^{*}v)^{T}]_{0}) = \xi_{i}v_{j} - \frac{\xi_{i}v_{i}}{2}\delta_{ij},
\ee
where $i,j$ are indices for $\dm$. For the mean curvature, one has $H_{\delta^{*}v}' = 
-\D f + v(H)$, so that the leading order term is just $-\D f$ with symbol $|\xi|^{2}f$. 
Hence one has elliptic data for the normal component $f$ of $v$. 

   For the tangential part of $v$, \eqref{3.5} gives 
$$\xi_{1}v_{2} = \xi_{2}v_{1} = 0 \ \ {\rm and} \ \ \xi_{1}v_{1} = \xi_{2}v_{2}.$$
Since $(\xi_{1}, \xi_{2}) \neq (0,0)$, it is elementary to see that the only solution 
of these equations is $v_{1} = v_{2} = 0$, which proves ellipticity. It is 
straightforward to verify further that $\Pi$ has Fredholm index 0.

{\endproof}

\bibliographystyle{plain}

\end{document}